\documentclass[a4paper]{article}

\usepackage{amsthm}
\usepackage{amssymb}
\usepackage{latexsym}
\usepackage{longtable}
\usepackage{epsfig}
\usepackage{amsmath}
\usepackage{hhline}

\usepackage{enumerate}

\newcommand{\sym}{\mathfrak{S}}

  \newcommand{\N}{\mathbb{N}}
  \newcommand{\Q}{\mathbb{Q}}
  \newcommand{\Z}{\mathbb{Z}}

\newtheorem{thm}{Theorem}[section]

\newtheorem{conj}[thm]{Conjecture}
\newtheorem{prop}[thm]{Proposition}
\long\def\symbolfootnote[#1]#2{\begingroup\def\thefootnote{\fnsymbol{footnote}}
\footnote[#1]{#2}\endgroup}

\begin{document}

\title{Defect of characters of the symmetric group}
\author{Jean-Baptiste Gramain\\ \\ \'Ecole Polytechnique F\'ed\'erale de Lausanne\\Lausanne,
Switzerland}
\date{September, 2008}
\maketitle

\begin{abstract}
Following the work of B. K\"ulshammer, J. B. Olsson and G. R. Robinson
on generalized blocks of the symmetric groups, we give a definition for the $\ell$-defect of characters of the symmetric group $\sym_n$, where $\ell >1$ is
an arbitrary integer. We prove that the $\ell$-defect is given by an analogue of the hook-length formula, and use it to prove, when $n < \ell^2$, an $\ell$-version of the McKay Conjecture in $\sym_n$.
\end{abstract}

\symbolfootnote[0]{2000 Mathematics Subject Classification 20C30 (primary), 20C15, 20C20 (secondary)}
\symbolfootnote[0]{EPFL-IGAT, B\^{a}timent de Chimie (BCH), CH-1015 Lausanne, Switzerland.}
\symbolfootnote[0]{email: jean-baptiste.gramain@epfl.ch}

\section{Introduction}

B. K\"{u}lshammer, J. B. Olsson and G. R. Robinson gave in \cite{KOR} a definition of {\emph{generalized blocks}} for a finite group. Let $G$ be a finite group, and denote by Irr$(G)$ the set of complex
irreducible characters of $G$. Take a union ${\cal C}$ of conjugacy classes of $G$ containing the identity. Suppose furthermore that ${\cal C}$ is {\emph{closed}}, that is, if $x \in {\cal C}$, and if $y \in G$ generates the same subgroup of $G$ as $x$, then $y \in {\cal C}$. For $\chi, \, \psi \in
{\mbox{Irr}}(G)$, we define the ${\cal C}${\emph{-contribution}} $\langle \chi, \, \psi\rangle _{\cal C}$ of $\chi$ and $\psi$ by
$$ \langle \chi, \, \psi \rangle_{\cal C}:= \displaystyle \frac{1}{|G|} \sum_{g \in {\cal C}} \chi(g) \psi(g^{-1}).$$
The fact that ${\cal C}$ is closed implies that, for any $\chi, \, \psi \in
{\mbox{Irr}}(G)$, $\langle \chi, \, \psi\rangle _{\cal C}$ is a rational number.

We say that $\chi, \, \psi \in {\mbox{Irr}}(G)$ belong to the same ${\cal C}${\emph{-block}} of $G$ if there exists a sequence of irreducible characters $\chi_1=\chi, \, \chi_2, \, \ldots , \, \chi_n= \psi$ of $G$ such that $\langle \chi_i, \, \chi_{i+1}\rangle _{\cal C} \neq 0$ for all $i \in \{1, \, \ldots , \, n-1 \}$. The ${\cal C}$-blocks define a partition of Irr$(G)$ (the fact that $1 \in {\cal C}$ ensures that each irreducible character of $G$ belongs to a ${\cal C}$-block). If we take
${\cal C}$ to be the set of $p$-regular elements of $G$ (i.e. whose
order is not divisible by $p$), for some prime $p$, then the ${\cal
C}$-blocks are just the ``ordinary" $p$-blocks (cf for example \cite{Navarro}, Theorem 3.19).

\medskip
Let $CF(G)$ be the set of complex class functions of $G$, and $\langle \, . \, , \, . \, \rangle$ be the ordinary scalar product on $CF(G)$. For any $\chi \in$ Irr$(G)$, we define $\chi^{\cal C} \in CF(G)$ by letting
$$\chi^{\cal C}(g)= \left\{ \begin{array}{ll} \chi(g) & \mbox{if} \; g \in {\cal C} \\ 0 & \mbox{otherwise} \end{array} \right. .$$
Then, for $\chi \in \mbox{Irr}(G)$, we have $\chi^{\cal C}= \sum_{\psi \in \mbox{Irr}(G)} \langle \chi^{\cal C}, \, \psi \rangle \psi=\sum_{\psi \in \mbox{Irr}(G)} \langle \chi, \, \psi \rangle_{\cal C} \psi$. Since $\langle \chi, \, \psi \rangle_{\cal C} \in \Q$ for all $\psi \in \mbox{Irr}(G)$, there exists $d \in \N$ such that $d \chi^{\cal C} $ is a generalized character of $G$. We call the smallest such positive integer the ${\cal C}${\emph{-defect}} of $\chi$, and denote it by $d_{\cal C} (\chi)$.

It is easy to check that $\chi \in \mbox{Irr}(G)$ has ${\cal C}$-defect 1 if and only if $\chi$ vanishes outside ${\cal C}$. This is also equivalent to the fact that $\{ \chi \}$ is a ${\cal C}$-block of $G$.

Writing $1_G$ for the trivial character of $G$, we see that, for any $\chi \in \mbox{Irr}(G)$, $d_{\cal C}(1_G) \chi^{\cal C}=\chi \otimes (d_{\cal C}(1_G) 1_G^{\cal C})$ is a generalized character, so that $d_{\cal C}(\chi)$ divides $d_{\cal C}(1_G)$. In particular, $1_G$ has maximal ${\cal C}$-defect.

Note that, if ${\cal C}$ is the set of $p$-regular elements of $G$ ($p$ a prime), then, for all $\chi \in$ Irr$(G)$, we have (cf for example \cite{Navarro}, Lemma 3.23) $d_{\cal C}(\chi)=\left(\frac{|G|}{\chi(1)}\right)_p=p^{d(\chi)}$, where $d(\chi)$ is the ordinary $p$-defect of $\chi$.

\medskip
One key notion defined in \cite{KOR} is that of {\emph{generalized
perfect isometry}}. Suppose $G$ and $H$ are finite
groups, and ${\cal C}$ and ${\cal D}$ are closed unions of conjugacy
classes of $G$ and $H$ respectively. Take $b$ a
union of ${\cal C}$-blocks of $G$, and $b'$ a
union of ${\cal D}$-blocks of $G$. A {\emph{generalized perfect isometry}}
between $b$ and $b'$ (with respect to ${\cal C}$ and ${\cal D}$) is a bijection with signs between $b$ and $b'$, which furthermore preserves contributions. That is, $I \colon b \longmapsto b'$ is a bijection such that, for each $\chi \in b$, there is a sign $\varepsilon (\chi )$, and such that
$$\forall \chi, \, \psi \in b, \, \; \langle I(\chi), \, I(\psi) \rangle _{{\cal D}}
=\langle \varepsilon (\chi) \chi, \, \varepsilon (\psi) \psi\rangle
_{{\cal C}}.$$
In particular, one sees that a generalized perfect isometry $I$ preserves the defect, that is, for all $\chi \in b$, we have $d_{\cal C}(\chi)=d_{\cal D}(I(\chi))$.

Note that, if ${\cal C}$ and ${\cal D}$ are the sets of $p$-regular elements of $G$ and $H$ respectively, then this notion is a bit weaker than that of {\emph{perfect isometry}} introduced by M. Brou\'e (cf \cite{Broue}). If two $p$-blocks $b$ and $b'$ are perfectly
isometric in Brou\'e's sense, then there is a generalized perfect
isometry (with respect to $p$-regular elements) between $b$ and $b'$.
It is however possible to exhibit generalized perfect isometries in
some cases where there is no perfect isometry in Brou\'e's sense (cf
\cite{JB}).

\medskip
K\"{u}lshammer, Olsson and Robinson defined and studied in \cite{KOR} the $\ell$-blocks of the
symmetric group, where $\ell \geq 2$ is any integer. They did this by taking ${\cal C}$ to be the set of $\ell${\emph{-regular}} elements, that is which have no cycle (in their canonical cycle decomposition) of length divisible by $\ell$
(in particular, if $\ell$ is a prime $p$, then the $\ell$-blocks are
just the $p$-blocks). 

\noindent
In section 2, we find the $\ell${\emph{-defect}} of the characters of the symmetric group $\sym_n$. It turns out (Theorem \ref{elldefectsn}) that it is given by an analogue of the hook-length formula (for the degree of a character). In section 3, we then use this to prove, when $n < \ell^2$, an $\ell$-analogue of the McKay Conjecture in $\sym_n$ (Theorem \ref{ellMcKay}).

\section{Hook-length formula}

\subsection{$\ell$-blocks of the symmetric group}

Take two integers $1 \leq \ell \leq n$, and consider the symmetric group $\sym_n$ on $n$ letters. The conjugacy classes and irreducible complex characters of $\sym_n$ are parametrized by the set $\{\lambda \vdash n \}$ of partitions of $n$. We write Irr$(\sym_n)=\{ \chi_{\lambda} , \, \lambda \vdash n\}$. An element of $\sym_n$ is said to be $\ell${\emph{-regular}} if none of its cycles has length divisible by $\ell$. We let ${\cal C}$ be the set of $\ell$-regular elements of $\sym_n$. The ${\cal C}$-blocks of $\sym_n$ are called $\ell${\emph{-blocks}}, and they satisfy the following

\begin{thm}[Generalized Nakayama Conjecture](\cite{KOR}, Theorem 5.13)\label{Nakayama}
Two characters $\chi_{\lambda}, \, \chi_{\mu} \in \mbox{Irr}(\sym_n)$ belong to the same $\ell$-block if and only if $\lambda$ and $\mu$ have the same $\ell$-core.
\end{thm}

The proof of this goes as follows. If $\langle \chi_{\lambda}, \, \chi_{\mu} \rangle \neq 0$, then an induction argument using the Murnaghan-Nakayama Rule shows that $\lambda$ and $\mu$ must have the same $\ell$-core. In particular, the partitions labeling the characters in an $\ell$-block all have the same $\ell$-weight, and we can talk about the $\ell${\emph{-weight}} of an $\ell$-block.

\smallskip

Conversely, let $B$ be the set of irreducible characters of $\sym_n$ labeled by those partitions of $n$ which have a given $\ell$-core, $\gamma$ say, and $\ell$-weight $w$. It is a well known combinatorial fact (cf for example \cite{James-Kerber}, Theorem 2.7.30) that the characters in $B$ are parametrized by the $\ell${\emph{-quotients}}, which can be seen as the set of  $\ell${\emph{-tuple of partitions}} of $w$. For $\chi_{\lambda} \in B$, the quotient $\beta_{\lambda}$ is a sequence $(\lambda^{(1)}, \, \ldots , \, \lambda^{(\ell )} )$ such that, for each $1 \leq i \leq \ell $, $\lambda^{(i)}$ is a partition of some $k_i$, $0 \leq k_i \leq w$, and $\sum_{i=1}^{\ell }k_i=w$ (the quotient $\beta_{\lambda}$ ``stores'' the information about how to remove $w$ $\ell$-hooks from $\lambda$ to get $\gamma$). We write $\beta_{\lambda} \Vdash w$. To prove that $B$ is an $\ell$-block of $\sym_n$, K\"{u}lshammer, Olsson and Robinson use a generalized perfect isometry between $B$ and the wreath product $\Z_{\ell} \wr \sym_w$ (where $\Z_{\ell}$ denotes a cyclic group of order $\ell$).

\smallskip

The conjugacy classes of $\Z_{\ell} \wr \sym_w$ are parametrized by the $\ell$-tuples of partitions of $w$ as follows (cf \cite{James-Kerber}, Theorem 4.2.8). Write $\Z_{\ell}= \{ g_1, \, \ldots , \, g_{\ell} \}$ the cyclic group of order $\ell$. The elements of the wreath product $\Z_{\ell} \wr \sym_w$ are
of the form $(h;\sigma)=(h_1, \, \ldots, \, h_w ; \sigma)$, with
$h_1, \, \ldots, \, h_w  \in \Z_{\ell} $ and $ \sigma \in \sym_w$. For any such
element, and for any $k$-cycle $\kappa=(j, j\kappa, \ldots, j
\kappa^{k-1})$ in $\sigma$, we define the {\emph{cycle product}} of $(h;\sigma)$
and $\kappa$ by
$$g((h;\sigma),\kappa)= h_j h_{j \kappa^{-1}} h_{j \kappa^{-2}} \ldots h_{j
\kappa^{-(k-1)}} \in \Z_{\ell}.$$If
$\sigma$ has cycle structure $\pi$ say, then we form $\ell$ partitions
$(\pi_1, \, \ldots, \, \pi_{\ell})$ from $\pi$ as follows: any cycle
$\kappa$ in $\pi$ gives a cycle of the same length in $\pi_i$ if $g((h;\sigma),\kappa)=g_i$. The
resulting $\ell$-tuple of partitions of $w$ describes the
{\emph{cycle structure}} of $(h;\sigma)$, and two elements of $\Z_{\ell} \wr \sym_w$ are conjugate if and only if they have
the same cycle structure. An element of $\Z_{\ell} \wr \sym_w$ is said to be {\emph{regular}} if it has no cycle product equal to 1.

\smallskip

The irreducible characters of $\Z_{\ell} \wr \sym_w$ are also canonically parametrized by the $\ell$-tuples of partitions of $w$ in the following way. Write Irr$(\Z_{\ell})= \{ \alpha_1, \, \ldots , \, \alpha_{\ell } \}$, and take $\beta_{\lambda}=(\lambda^{(1)},  \, \ldots , \, \lambda^{(\ell )} ) \Vdash w$, with $\lambda^{(i)} \vdash k_i$ as above ($1 \leq i \leq \ell$). The irreducible character $\alpha_1^{k_1} \otimes \cdots  \otimes \alpha_{\ell}^{k_{\ell}}$ of the {\emph{base group}} $\Z_{\ell}^w$ can be extended in a natural way to its inertia subgroup $(\Z_{\ell} \wr \sym_{k_1}) \times \cdots \times (\Z_{\ell} \wr \sym_{k_{\ell}})$, giving the irreducible character $\prod_{i=1}^{\ell} \widehat{ \alpha_i^{k_i} }$. The tensor product $\prod_{i=1}^{\ell} \widehat{ \alpha_i^{k_i} } \otimes \chi_{\lambda^{(i)}}$ is an irreducible character of $(\Z_{\ell} \wr \sym_{k_1}) \times \cdots \times (\Z_{\ell} \wr \sym_{k_{\ell}})$ which extends $\prod_{i=1}^{\ell} \widehat{ \alpha_i^{k_i} }$, and it remains irreducible when induced to $\Z_{\ell} \wr \sym_w$. We denote by $\chi_{\beta_{\lambda}}$ this induced character. Furthermore, any irreducible character of $\Z_{\ell} \wr \sym_w$ can be obtained in this way.

\smallskip

In \cite{KOR}, the authors show that the map $\chi_{\lambda} \longmapsto \chi_{\beta_{\lambda}}$ is a generalized perfect isometry between $B$ and Irr$(\Z_{\ell} \wr \sym_w)$, with respect to $\ell$-regular elements of $\sym_n$ and regular elements of $\Z_{\ell} \wr \sym_w$.

On the other hand, they show that, writing $reg$ for the set of regular elements of $\Z_{\ell} \wr \sym_w$, we have, for all $\chi \in \mbox{Irr}(\Z_{\ell} \wr \sym_w)$,
\begin{equation} \label{def} \Z \ni \displaystyle \frac{\ell ^w w! \langle \chi, \, 1_{\Z_{\ell} \wr \sym_w} \rangle_{reg}}{\chi (1)} \equiv (-1)^w \; (\mbox{mod} \, \ell),\end{equation}
where $1_{\Z_{\ell} \wr \sym_w}$ is the trivial character of $\Z_{\ell} \wr \sym_w$. In particular, $ \langle \chi, \, 1_{\Z_{\ell} \wr \sym_w} \rangle_{reg} \neq 0$. Using the generalized perfect isometry we described above, this implies that there exists a character $\chi_{\lambda} \in B$ such that, for all $\chi_{\mu} \in B$, we have $\langle \chi_{\lambda} , \, \chi_{\mu} \rangle_{\cal C} \neq 0$, where ${\cal C}$ is the set of $\ell$-regular elements of $\sym_n$. In particular, all the characters in $B$ belong to the same $\ell$-block of $\sym_n$, which ends the proof of Theorem 2.1.

\subsection{$\ell$-defect of characters}

Using the ingredients in the proof of Theorem 2.1, we can now compute explicitly the $\ell${\emph{-defects}} of the irreducible characters of $\sym_n$ (that is, their ${\cal C}$-defect, where ${\cal C}$ is the set of $\ell$-regular elements of $\sym_n$).

As we remarked earlier, if $\lambda$ is a partition of $n$ of $\ell$-weight $w$, then, because of the generalized perfect isometry we described above, the $\ell$-defect $d_{\ell}(\chi_{\lambda})$ of $\chi_{\lambda} \in \mbox{Irr}(\sym_n)$ is the same as the $reg$-defect $d_{reg}(\chi_{\beta_{\lambda}})$ of $\chi_{\beta_{\lambda}} \in \mbox{Irr}(\Z_{\ell} \wr \sym_w)$, where $\beta_{\lambda}$ is the $\ell$-quotient of $\lambda$. It is in fact these $reg$-defects we will compute.

First note that, if $w=0$, then $\lambda$ is its own $\ell$-core, so that $\chi_{\lambda}$ is alone in its $\ell$-block, and $d_{\ell}(\chi_{\lambda})=1$. We therefore now fix $w \geq 1$.

\medskip
We write $\pi$ the set of primes dividing $\ell$. Every positive integer $m$ can be factorized uniquely as $m=m_{\pi}m_{\pi'}$, where every prime factor of $m_{\pi}$ belongs to $\pi$ and no prime factor of $m_{\pi'}$ is contained in $\pi$. We call $m_{\pi}$ the $\pi$-part of $m$. 

Using results of Donkin (cf \cite{Donkin}) and equality (\ref{def}), K\"{u}lshammer, Olsson and Robinson proved the following

\begin{thm}(\cite{KOR}, Theorem 6.2)
The $reg$-defect of the trivial character of $\Z_{\ell} \wr \sym_w$ is $\ell ^w w!_{\pi}$.
\end{thm}

In particular, since $1_{\Z_{\ell} \wr \sym_w}$ has maximal $reg$-defect, we see that, for any $\chi \in \mbox{Irr}(\Z_{\ell} \wr \sym_w)$, $d_{reg}(\chi)$ is a $\pi$-number.

\medskip
We can now compute the $reg$-defect of any irreducible character $\chi$ of $\Z_{\ell} \wr \sym_w$. It turns out that it is sufficient to know the $reg$-contribution of $\chi$ with the trivial character, and this is given by (\ref{def}). We have the following

\begin{prop}
Take any integers $ \ell \geq 2$ and $w \geq 1$. Then, for any \\$\chi \in \mbox{Irr}(\Z_{\ell} \wr \sym_w)$, we have $d_{reg}(\chi)= \displaystyle \frac{\ell^w (w!)_{\pi}}{\chi(1)_{\pi}}$.
\end{prop}

\begin{proof}
Take $\chi \in \mbox{Irr}(\Z_{\ell} \wr \sym_w)$. Recall that, by (\ref{def}),
$$\Z \ni \displaystyle \frac{\ell ^w w! }{\chi (1)}  \langle\chi,1\rangle_{reg} \equiv (-1)^w \; (\mbox{mod} \, \ell) .$$
Now $d_{reg}(\chi)$ is a $\pi$-number, so that  $\langle\chi, \, 1\rangle_{reg}$ is a rational whose (reduced) denominator is a $\pi$-number. This implies that $ \frac{\ell ^w (w!)_{\pi} }{\chi (1)_{\pi}}  \langle\chi, \, 1\rangle_{reg} \in \Z$.

Furthermore, from (\ref{def}), we also deduce that, {\bf for each $p \in \pi$},
$$ \displaystyle \frac{\ell ^w w! }{\chi (1)}  \langle\chi, \, 1\rangle_{reg} \not \equiv 0 \; (\mbox{mod} \, p).$$
Thus, for any $p \in \pi$, $ \displaystyle \frac{\ell ^w (w!)_{\pi} }{\chi (1)_{\pi}}  \langle\chi, \, 1\rangle_{reg} \not \equiv 0 \; (\mbox{mod} \, p)$. Hence $ \frac{\ell ^w (w!)_{\pi} }{\chi (1)_{\pi}} $ is the {\bf smallest} positive integer $d$ such that $d \langle\chi, \, 1\rangle_{reg}  \in \Z$. This implies that $ \frac{\ell ^w (w!)_{\pi} }{\chi (1)_{\pi}}$ divides $d_{reg}(\chi)$ (indeed, by definition, $d_{reg}(\chi) \langle\chi, \, 1\rangle_{reg}  \in \Z$, and $d_{reg}(\chi)$ is a $\pi$-number).

\smallskip
Now, conversely, if $\psi \in \mbox{Irr}(\Z_{\ell} \wr \sym_w)$, then $\langle\chi, \, \psi\rangle_{reg} \in \Q$, so also (since $\chi(1)$ divides $|\Z_{\ell} \wr \sym_w|=\ell^w w!$) $ \frac{\ell ^w w! }{\chi (1)}  \langle\chi, \, \psi\rangle_{reg} \in \Q$. However,
$$\displaystyle \frac{\ell ^w w! }{\chi (1)}  \langle\chi, \, \psi\rangle_{reg}= \frac{\ell^w w!}{\ell^w w!} \sum_{g \in reg/ \sim} \frac{K_g \chi(g)}{\chi(1)} \psi(g^{-1})$$
(where the sum is taken over representatives for the regular classes, and, for $g$ such a representative, $K_g$ is the size of the conjugacy class of $g$). And, for each $g$ in the sum,  $\frac{K_g \chi(g)}{\chi(1)}$ and $ \psi(g^{-1})$ are both algebraic integers. Hence $ \frac{\ell ^w w! }{\chi (1)}  \langle\chi, \, \psi\rangle_{reg}$ is also an algebraic integer, and thus an integer. Thus
$$ \forall \psi \in \mbox{Irr}(\Z_{\ell} \wr \sym_w) , \; \; \displaystyle \frac{\ell ^w w! }{\chi (1)}  \langle\chi, \, \psi\rangle_{reg} \in \Z$$
and this implies that $d_{reg}(\chi)$ divides $\frac{\ell ^w w! }{\chi (1)} $, and, $d_{reg}(\chi)$ being a $\pi$-number, $d_{reg}(\chi)$ divides $ \frac{\ell ^w (w!)_{\pi} }{\chi (1)_{\pi}}$. Hence we finally get $d_{reg}(\chi)=  \frac{\ell ^w (w!)_{\pi} }{\chi (1)_{\pi}}$.
\end{proof}

We want to express the $\ell$-defect of a character in terms of hook-lengths. For any $\lambda \vdash n$, we write ${\cal H}(\lambda)$ for the set of hooks in $\lambda$, and ${\cal H}_{\ell}(\lambda)$ for the set of hooks in $\lambda$ whose length is divisible by $\ell$. Similarly, if $\beta_{\lambda} = (\lambda^{(1)}, \, \ldots , \, \lambda^{(\ell)}) \Vdash w$, we define a {\emph{hook}} in $\beta_{\lambda}$ to be a hook in any of the $\lambda^{(i)}$'s, and write ${\cal H}(\beta_{\lambda})$ for the set of hooks in $\beta_{\lambda}$. Finally, for any hook $h$ (in a partition or a tuple of partitions), we write $|h|$ for the length of $h$.

We will use the following classical results about hooks (cf for example \cite {James-Kerber}, \S 2.3 and \S2.7)

\begin{thm}\label{HLF}
Let $n \geq \ell \geq 2$ be any two integers, and let $\lambda$ be any partition of $n$. Then

\begin{enumerate}[(i)]

\item 
(Hook-Length Formula, \cite {James-Kerber}, Theorem 2.3.21) We have $ \frac{|\sym_n|}{\chi_{\lambda}(1)} = \prod_{h \in {\cal H}(\lambda)} |h|$,

\item (\cite {James-Kerber}, 2.7.40) if $\lambda$ has $\ell$-weight $w$, then $|{\cal H}_{\ell}(\lambda)|=w$,

\item (\cite {James-Kerber}, Lemma 2.7.13 and Theorem 2.7.16) if $\beta_{\lambda}$ is the $\ell$-quotient of $\lambda$, then $\{ |h|, \, h \in {\cal H}_{\ell}(\lambda) \} = \{ \ell |h'|, \, h' \in {\cal H}(\beta_{\lambda}) \}$.

\end{enumerate}

\end{thm}

We can now establish the following

\begin{prop} 
If $n \geq \ell \geq 2$ are any two integers, $\pi$ is the set of primes dividing $\ell$, and $\lambda \vdash n$ has $\ell$-weight $w \neq 0$ and $\ell$-quotient $\beta_{\lambda}$, then
$$\displaystyle \frac{\ell^w (w!)_{\pi}}{\chi_{\beta_{\lambda}}(1)_{\pi}} = \prod_{h \in {\cal H}_{\ell}(\lambda)} |h|_{\pi}.$$
\end{prop}

\begin{proof}
Write $\beta_{\lambda} = (\lambda^{(1)}, \, \ldots , \, \lambda^{(\ell)})$, where $\lambda^{(i)} \vdash k_i$ for $1 \leq i \leq \ell$. First note that, by construction of $\chi_{\beta_{\lambda}}$, and since the irreducible characters of $\Z_{\ell}$ all have degree $1$, we have
$$\chi_{\beta_{\lambda}}(1)=\displaystyle \frac{\ell^w w!}{\prod_{i=1}^{\ell} \ell^{k_i} k_i !} \chi_{\lambda^{(1)}}(1) \ldots \chi_{\lambda^{(\ell)}}(1).$$
Thus, by the Hook-Length Formula (Theorem \ref{HLF}, (i)),
$$\chi_{\beta_{\lambda}}(1)=\displaystyle \frac{w!}{\prod_{h \in {\cal H}(\beta_{\lambda})} |h|} \; \; \mbox {and} \; \;  \frac{|\Z_{\ell} \wr \sym_w |}{\chi_{\beta_{\lambda}}(1)}= \ell^w \prod_{h \in {\cal H}(\beta_{\lambda})} |h|.$$
We therefore get
$$\displaystyle \frac{\ell^w (w!)_{\pi}}{\chi_{\beta_{\lambda}}(1)_{\pi}} =  \frac{|\Z_{\ell} \wr \sym_w |_{\pi}}{\chi_{\beta_{\lambda}}(1)_{\pi}}= \ell^w \prod_{h \in {\cal H}(\beta_{\lambda})} |h|_{\pi}.$$
Now, by Theorem \ref{HLF} (ii) and (iii), we have $|{\cal H}(\beta_{\lambda})| =w$, so that $\ell^w  \prod_{h \in {\cal H}(\beta_{\lambda})} |h| = \prod_{h \in {\cal H}(\beta_{\lambda})} \ell |h|$, and, by Theorem 2.4 (iii), $  \prod_{h \in {\cal H}(\beta_{\lambda})} \ell |h| =  \prod_{h \in {\cal H}_{\ell}(\lambda)}  |h|$. Taking $\pi$-parts, we obtain $ \frac{\ell^w (w!)_{\pi}}{\chi_{\beta_{\lambda}}(1)_{\pi}} = \prod_{h \in {\cal H}_{\ell}(\lambda)} |h|_{\pi}$, as announced.
\end{proof}

Combining Propositions 2.3 and 2.5, we finally get

\begin{thm}\label{elldefectsn}
Let $2 \leq \ell \leq n$ be any two integers, and let $B$ be an $\ell$-block of $\sym_n$ of weight $w$. Then
\begin{enumerate}[(i)]
\item If $w=0$, then $B= \{ \chi_{\lambda} \}$ for some partition $\lambda$ of $n$, and $d_{\ell}(\chi_{\lambda})=1$.

\item If $w > 0$, and if $\chi_{\lambda} \in B$, then $d_{\ell}(\chi_{\lambda})= \prod_{h \in {\cal H}_{\ell}(\lambda)} |h|_{\pi}$, where $\pi$ is the set of primes dividing $\ell$ (that is, $d_{\ell}(\chi_{\lambda})$ is the $\pi$-part of the product of the hook-lengths divisible by $\ell$ in $\lambda$).

\end{enumerate}
\end{thm}

\section{McKay Conjecture}

\subsection{McKay Conjecture, generalization}
In this section, we want to study an $\ell$-analogue of the following

\begin{conj}[McKay]
Let $G$ be a finite group, $p$ be a prime, and $P$ be a Sylow $p$-subgroup of $G$. Then the numbers of irreducible complex characters whose degree is not divisible by $p$ are the same for $G$ and $N_G(P)$.
\end{conj}

First note that the McKay Conjecture was proved by J. B. Olsson for the symmetric group (cf \cite{Olsson}). In order to generalize this to an arbitrary integer $\ell$, we will use the results of \cite{JB2}, which we summarize here. Let $2 \leq \ell \leq n$ be integers. Suppose furthermore that $n < \ell^2$, and write $n = \ell w + r$, with $0 \leq w, \, r < \ell$. We define a {\emph{Sylow $\ell$-subgroup}} of $\sym_n$ to be any subgroup of $\sym_n$ generated by $w$ disjoint $\ell$-cycles. In particular, if $\ell$ is a prime $p$, then the Sylow $\ell$-subgroups of $\sym_n$ are just its Sylow $p$-subgroups. Then any two Sylow $\ell$-subgroups of $\sym_n$ are conjugate, and they are Abelian. Let ${\cal L}$ be a Sylow $\ell$-subgroup of $\sym_n$. In \cite{JB2}, a notion of $\ell${\emph{-regular}} element is given, which coincide with the notion of $p$-regular element if $\ell$ is a prime $p$. Using this, one can construct the $\ell${\emph{-blocks}} of $N_{\sym_n}({\cal L})$, and show that they satisfy an analogue of Brou\'e's Abelian Defect Conjecture (cf \cite{JB2}, Theorem 4.1). We will show that, still in the case where $n < \ell^2$, an analogue of the McKay Conjecture also holds. However, if we just replace $p$ by any integer $\ell$, and consider irreducible characters of degree not divisible by $\ell$, or even coprime to $\ell$, then the numbers differ in $\sym_n$ and $N_{\sym_n}({\cal L})$. Instead, we will use the notion of $\ell$-defect, and prove that the numbers of irreducible characters of maximal $\ell$-defect are the same in $\sym_n$ and $N_{\sym_n}({\cal L})$ (note that, if $\ell$ is a prime, then both statements coincide).

\subsection{Defect and weight}

In order to study characters of $\sym_n$ of maximal $\ell$-defect, we need the following result, which tells us where to look for them.

\begin{prop}\label{defectmax}
Let $\ell \geq 2$ and $0 \leq w, \, r < \ell$ be any integers, and let $\lambda$ be a partition of $n= \ell w + r$. If $\chi_{\lambda} \in \mbox{Irr}(\sym_n)$ has maximal $\ell$-defect, then $\lambda$ has (maximal) $\ell$-weight $w$.
\end{prop}

\begin{proof}
First note that, if $\ell$ is a prime, then this can be proved in a purely arithmetic way (cf \cite{MacDonald}). This doesn't seem to be the case when $\ell$ is no longer a prime, and we will use the abacus instead. For a complete description of the abacus, we refer to \cite{James-Kerber}, \S 2.7 (note however that the abacus we use here is the horizontal mirror image of that described by James and Kerber).

Suppose, for a contradiction, that $\lambda$ has $\ell$-weight $v < w$. By the previous section, the $\ell$-defect of $\chi_{\lambda}$ is the $\pi$-part of the product of the hook-lengths divisible by $\ell$ in $\lambda$. Now these are visible on the $\ell$-abacus of $\lambda$. This has $\ell$ runners, and a hook of length $k \ell$ ($k \geq 1$) corresponds to a bead sitting, on a runner, $k$ places above an empty spot. In particular, the $(\ell)$-hooks (i.e. those whose length is divisible by $\ell$) in $\lambda$ are stored on at most $v$ runners. To establish the result, we will construct a partition $\mu$ of $n$ of weight $w$, and such that $d_{\ell}(\chi_{\mu}) > d_{\ell}(\chi_{\lambda})$.

Start with the $\ell$-abacus of any partition $\nu$ of $r$. On the (at most) $v$ runners used by $\lambda$, take some beads up to encode the same $(\ell)$-hooks as for $\lambda$. Then, on $w-v$ of the (at least) $\ell-v > w -v$ remaining runners, take the highest bead one place up. The resulting abacus then corresponds to a partition of $n = r + \ell w= r + \ell v + \ell (w-v)$, and we see that $d_{\ell}(\chi_{\mu}) = \ell^{w-v} d_{\ell}(\chi_{\lambda})$ (indeed, the $(\ell)$-hooks in $\mu$ are precisely those in $\lambda$, together with $w-v$ hooks of length $\ell$). This proves the result.

\end{proof}

\subsection{Generalized Perfect Isometry}

We describe here the analogue of Brou\'e's Abelian Defect Conjecture given in \cite{JB2} (Theorem 4.1). We take any integers $\ell \geq 2$ and $0 \leq w, \, r <\ell$, and $G=\sym_{\ell w + r}$. We take an Abelian {\emph{Sylow $\ell$-subgroup}} ${\cal L}$ of $G$; that is, ${\cal L} \cong \Z_{\ell}^w$ is generated by $w$ disjoint $\ell$-cycles. Then ${\cal L}$ is a natural subgroup of $\sym_{\ell w}$, and we have $N_G({\cal L}) \cong N_{\sym_{\ell w}}({\cal L}) \times \sym_r$ and Irr($N_G({\cal L})) = \mbox{Irr} (N_{\sym_{\ell w}}({\cal L})) \otimes \mbox{Irr}(\sym_r)$. Now $N_{\sym_{\ell w}}({\cal L}) \cong N \wr \sym_w = N_{\sym_{\ell }}(L) \wr \sym_w$, where $L=\langle \pi \rangle \cong \Z_{\ell}$ is (a subgroup of $\sym_{\ell}$) generated by a single $\ell$-cycle. As in the sketch of the proof of Theorem \ref{Nakayama}, we see that the conjugacy classes and irreducible characters of $N_{\sym_{\ell w}}({\cal L})$ are parametrized by the $s$-tuples of partitions of $w$, where $s$ is the number of conjugacy classes of $N$. Among these, there is a unique conjugacy class of $\ell$-cycles, for which we take representative $\pi$. We take representatives $\{ g_1=\pi, \, g_2, \, \ldots , \, g_s\}$ for the conjugacy classes of $N$. Considering as $\ell${\emph{-regular}} any element of $N$ not conjugate to the $\ell$-cycle $\pi$, we can construct the $\ell${\emph{-blocks}} of $N$, and show that the principal $\ell$-block contains $\ell$ characters, which we label $\psi_1, \, \ldots , \, \psi_{\ell}$, and that each of the remaining $s-\ell$ characters, labeled $\psi_{\ell + 1}, \, \ldots , \, \psi_s$, is alone in its $\ell$-block (cf \cite{JB2}, Section 2). Using the construction presented after Theorem \ref{Nakayama}, we label the conjugacy classes and irreducible characters of $N \wr \sym_w$ by the $s$-tuples of partitions of $w$. An element of $N \wr \sym_w$ of cycle type $(\pi_1, \, \ldots , \, \pi_s) \Vdash w$ is called $\ell${\emph{-regular}} if $\pi_1=\emptyset$ (and $\ell${\emph{-singular}} otherwise). Then one shows that the $\ell$-blocks of $N \wr \sym_w$ are the principal $\ell$-block, $b_0= \{ \chi^{\alpha}, \, \alpha=(\alpha_1, \, \ldots , \, \alpha_{\ell}, \, \emptyset, \, \ldots, \, \emptyset \} \Vdash w\}$, and blocks of size $1$, $\{ \chi^{\alpha} \}$, whenever $\alpha \Vdash w$ is such that $\alpha_k \neq \emptyset$ for some $\ell < k \leq s$ (\cite{JB2}, Theorem 3.7 and Corollary 3.11).

Finally, an element of $N_G({\cal L}) \cong N_{\sym_{\ell w}}({\cal L}) \times \sym_r$ is said to be $\ell${\emph{-regular}} if its $ N_{\sym_{\ell w}}({\cal L})$-part is $\ell$-regular in the above sense (so that, if $\ell$ is a prime $p$, then the notions of $\ell$-regular and $p$-regular coincide). Then, we can summarize the results of \cite{JB2} as follows:

\begin{thm}[\cite{JB2}, Theorem 4.1]\label{genperfisom}
Let the notation be as above. Then any $\ell$-block of $N_G({\cal L})$ has size $1$ or belongs to $\{ b_0 \otimes \{ \psi \} , \, \psi \in \mbox{Irr}(\sym_r) \}$. Furthermore, for any $ \psi \in \mbox{Irr}(\sym_r) $, there is a generalized perfect isometry (with respect to $\ell$-regular elements) between $b_0 \otimes \{ \psi \}$ and $B_{\psi}$, where $B_{\psi}$ is the $\ell$-block of $\sym_{\ell w +r}$ consisting of the irreducible characters labeled by partitions with $\ell$-core $\psi$.
\end{thm}
\noindent
(Note that any partition of $r$ does appear as $\ell$-core of a partition of $\ell w + r$ (for example, if $\gamma \vdash r$, then $\gamma$ is the $\ell$-core of $(\gamma, \, 1^{\ell w}) \vdash \ell w + r$).)

\subsection{Analogues of the McKay Conjecture}

We can now give the analogue of the McKay Conjecture we announced. Let, as before, $\ell \geq 2$ and $0 \leq w, \, r <\ell$ be integers, $n=\ell w + r$, and ${\cal L}$ be an Abelian Sylow $\ell$-subgroup of $\sym_n$. By Proposition \ref{defectmax}, any irreducible character of $\sym_n$ of maximal $\ell$-defect has (maximal) $\ell$-weight $w$, hence belongs to one of the $B_{\psi}$'s, $\psi \in \mbox{Irr}(\sym_r)$. Since any generalized perfect isometry preserves the defect, Theorem \ref{genperfisom} provides a bijection between the sets of irreducible characters of maximal $\ell$-defect and $\ell$-weight $w$ of $\sym_n$ and of characters of maximal $\ell$-defect in $N_{\sym_n}({\cal L})$. We therefore obtain

\begin{thm}\label{ellMcKay}
With the above notations, the numbers of irreducible characters of maximal $\ell$-defect are the same in $\sym_n$ and $N_{\sym_n}({\cal L})$.
\end{thm}
\noindent
{\bf{Remark:}} furthermore, we have an explicit bijection, essentially given by taking $\ell$-quotients of partitions.

\medskip
\noindent
In fact, Theorem \ref{genperfisom} gives something a bit stronger, namely
\begin{thm}For any $\ell$-defect $\delta \neq 1$, there is a bijection between the set of irreducible characters of $\sym_n$ of $\ell$-weight $w$ and $\ell$-defect $\delta$ and the set of irreducible characters of $N_{\sym_n}({\cal L})$ of $\ell$-defect $\delta$.
\end{thm}

Now, McKay's Conjecture is stated (and, in the case of symmetric groups, proved) without any hypothesis on the Sylow $p$-subgroups. One would therefore want to generalize the above results to the case where $n \geq \ell^2$. Examples seem to indicate that such analogues do indeed hold in this case, and that a bijection is given by taking, not only the $\ell$-quotient of a partition, but its $\ell${\emph{-tower}} (cf \cite{Olsson}).

In order to prove these results, one would first need to generalize Proposition \ref{defectmax}, showing that, for any $2 \leq \ell \leq n$, if $\chi_{\lambda} \in \mbox{Irr}(\sym_n)$ has maximal $\ell$-defect, then $\lambda$ has maximal $\ell$-weight, but also maximal $\ell^2$-weight, maximal $\ell^3$-weight, ... . If $\ell$ is a prime, then this is known to be true (cf \cite{MacDonald}). However, it seems hard to prove in general, even when $n= \ell^2$. The particular case where $\ell$ is squarefree is much easier.

Also, one would need to generalize the results of \cite{JB2}, while making sure that, when $\ell$ is a prime $p$, the notions of $\ell$-regular and $p$-regular elements still coincide.

\bigskip

\bigskip
\emph{Jean-Baptiste Gramain}

\emph{EPFL}

\emph{IGAT}

\emph{B\^{a}timent de Chimie (BCH)}

\emph{CH-1015 Lausanne}

\emph{Switzerland}

\medskip
jean-baptiste.gramain@epfl.ch

\end{document}